\documentclass{article}
\usepackage{amsmath,amssymb}
\newtheorem{df}{Definition}[section]
\newtheorem{thm}[df]{Theorem}

\newtheorem{cor}[df]{Corollary}
\title{Continuous wavelet transforms\\ on $n$--dimensional spheres}
\author{I. Iglewska-Nowak\footnote{West Pomeranian University of Technology in Szczecin, School of Mathematics, al. Piast\'ow 17, 70--310 Szczecin, Poland}}

\begin{document}

\maketitle

\bibliographystyle{amsplain}

\begin{abstract} In this paper, we are concerned with $n$--dimensional spherical wavelets derived from the theory of approximate identities. For nonzonal bilinear wavelets introduced by Ebert \emph{et al.} in 2009 we prove isometry and Euclidean limit property. Further, we develop a theory of linear wavelets. In the end, we discuss the relationship to other wavelet constructions.
\end{abstract}

\begin{bfseries}Keywords:\end{bfseries} spherical wavelets, $n$-spheres, singular integrals, approximate identities \\
\begin{bfseries}AMS Classification: 42C40, 42B20\end{bfseries}

\section{Introduction}

Continuous wavelet transform based on singular integrals on the $2$--sphere was introduced 1996 in~\cite{FW} (compare also \cite{FW-C} and \cite{FGS-book}). The wavelets are zonal functions (they are in principle derivatives of kernels of singular integrals) such that the parameter space of the wavelet transform is $\mathcal S^2\times\mathbb R_+$. Both linear and bilinear wavelet transforms are considered, where \emph{linear} means that the decomposition is done by the wavelet, and no wavelet is needed for the reconstruction, whereas in the \emph{bilinear} case one uses a wavelet both for the analysis and for the synthesis, compare also \cite[Remark on p. 232]{FGS-book}.

A generalization to three and more dimensions was done by Bernstein \emph{at al.} starting from 2009, cf. \cite{sB09}, \cite{BE10WS3}, and \cite{BE10KBW}. A special case of diffusive wavelets (i.e., satisfying an additional condition of diffusivity) is studied in Ebert's Ph.D.--thesis~\cite{sE11}.

In \cite{BE10WS3} a bilinear wavelet transform over $SO(3)$ is introduced with wavelets zonal with respect to an $SO(3)$ transform $g_0$. This idea is further developed in~\cite{EBCK09} to nonzonal bilinear wavelets over $SO(n)$. In this paper, we prove some properties of the nonzonal bilinear wavelet transform over $SO(n)$ and we introduce a nonzonal linear wavelet transform over $SO(n)$. Finally, we discuss the relationship between these and other wavelet constructions.

The paper is organized as follows. Section~\ref{sec:sphere} contains basic information about analysis of functions on spheres. The notions of singular integrals and approximate identities are recapitulated and illustrated on some examples. In Subsection~\ref{subs:dilations} we show that dilations via stereographic projections that were used in~\cite{AVn} create a kernel of an approximate identity. Section~\ref{sec:bilinearwv} digests and completes information about nonzonal and zonal bilinear wavelet transform introduced in~\cite{EBCK09}, especially it is shown that spherical wavelets under some mild additional conditions have the Euclidean limit property, i.e., for small scales they behave like wavelets over Euclidean space. In Section~\ref{sec:linearwv} we introduce linear spherical wavelets and discuss their properties. Relationship between the presented constructions and other wavelets used for analysis of spherical functions is briefly discussed in Section~\ref{sec:other}.

\section{Preliminaries}\label{sec:sphere}

\subsection{Functions on the sphere}

By $\mathcal{S}^n$ we denote the $n$--dimensional unit sphere in $n+1$--dimensional Euclidean space~$\mathbb{R}^{n+1}$ with the rotation--invariant measure~$d\sigma$ normalized such that
$$
\Sigma_n=\int_{\mathcal{S}^n}d\sigma=\frac{2\pi^{(n+1)/2}}{\Gamma\bigl((n+1)/2\bigr)}.
$$
The surface element $d\sigma$ is explicitly given by
$$
d\sigma=\sin^{n-1}\theta_1\,\sin^{n-2}\theta_2\dots\sin\theta_{n-1}d\theta_1\,d\theta_2\dots d\theta_{n-1}d\varphi,
$$
where $(\theta_1,\theta_2,\dots,\theta_{n-1},\varphi)\in[0,\pi]^{n-1}\times[0,2\pi)$ are spherical coordinates satisfying
\begin{align*}
x_1&=\cos\theta_1,\\
x_2&=\sin\theta_1\cos\theta_2,\\
x_3&=\sin\theta_1\sin\theta_2\cos\theta_3,\\
&\dots\\
x_{n-1}&=\sin\theta_1\sin\theta_2\dots\sin\theta_{n-2}\cos\theta_{n-1},\\
x_n&=\sin\theta_1\sin\theta_2\dots\sin\theta_{n-2}\sin\theta_{n-1}\cos\varphi,\\
x_{n+1}&=\sin\theta_1\sin\theta_2\dots\sin\theta_{n-2}\sin\theta_{n-1}\sin\varphi.
\end{align*}
$\left<x,y\right>$ or $x\cdot y$ stands for the scalar product of vectors with origin in~$O$ and endpoint on the sphere. As long as it does not lead to misunderstandings, we identify these vectors with points on the sphere.

By $\mathcal{X}(\mathcal{S}^n)$ we denote the space $\mathcal{C}(\mathcal{S}^n)$ or $\mathcal{L}^p(\mathcal{S}^n)$, \mbox{$1\leq p<\infty$}, with norm given by
$$
\|f\|_{\mathcal{C}(\mathcal{S}^n)}=\sup_{x\in\mathcal S^n}|f(x)|
$$
or
$$
\|f\|_{\mathcal{L}^p(\mathcal{S}^n)}=\left[\frac{1}{\Sigma_n}\int_{\mathcal S^n}|f(x)|^p\,d\sigma(x)\right]^{1/p},
$$
respectively, and by $\mathcal{X}\bigl([-1,1]\bigr)$ --- the space $\mathcal{C}\bigl([-1,1]\bigr)$ or $\mathcal{L}_\lambda^p\bigl([-1,1]\bigr)$, \mbox{$1\leq p<\infty$},
 with norm given by
$$
\|f\|_{\mathcal{C}([-1,1])}=\sup_{t\in[-1,1]}|f(t)|
$$
or
$$
\|f\|_{\mathcal{L}_\lambda^p([-1,1])}=\left[\frac{\Sigma_{2\lambda}}{\Sigma_{2\lambda+1}}\int_{-1}^1|f(t)|^p\left(1-t^2\right)^{\lambda-1/2}dt\right]^{1/p},
$$
respectively. The scalar product of $f,g\in\mathcal L^2(\mathcal S^n)$ is defined by
$$
\left<f,g\right>_{\mathcal L^2(\mathcal S^n)}=\frac{1}{\Sigma_n}\int_{\mathcal S^n}\overline{f(x)}\,g(x)\,d\sigma(x),
$$
such that $\|f\|_2^2=\left<f,f\right>$.

A function is called zonal if its value depends only on $\theta=\theta_1=\left<\hat e,x\right>$, where~$\hat e$ is the north pole of the sphere
$$
\hat e=(1,0,0,\dots,0).
$$
It is invariant with respect to the rotation about the axis through~$O$ and~$\hat e$. The subspace of $p$--integrable zonal functions is isomorphic to and will be identified with the space~$\mathcal L_\lambda^p$, and the norms satisfy
$$
\|f\|_{\mathcal{L}^p(\mathcal{S}^n)}=\|f\|_{\mathcal{L}^p([-1,1])},
$$
where~$\lambda$ and~$n$ are related by
$$
\lambda=\frac{n-1}{2}.
$$
We identify zonal functions with functions over the interval $[-1,1]$, i.e., whenever it does not lead to mistakes, we write
$$
f(x)=f(\cos\theta_1).
$$

Gegenbauer polynomials $C_l^\lambda$ of order~$\lambda\in\mathbb R$ and degree $l\in\mathbb{N}_0$ are defined in terms of their generating function
$$
\sum_{l=0}^\infty C_l^\lambda(t)\,r^l=\frac{1}{(1-2tr+r^2)^\lambda},\qquad t\in[-1,1],
$$
they also satisfy the relation
$$
\sum_{l=0}^\infty\frac{\lambda+l}{\lambda}\,C_l^\lambda(t)\,r^l=\frac{1-r^2}{(1-2tr+r^2)^{\lambda+1}},\qquad t\in[-1,1].
$$
(Legendre polynomials are Gegenbauer polynomials with $\lambda=1/2$.) Explicitly we have
\begin{equation}\label{eq:C_jako_szereg}
C_l^\lambda(t)=\sum_{j=0}^{[l/2]}(-1)^j\frac{\Gamma(\lambda+l-j)}{j!\,(l-2j)!\,\Gamma(\lambda)}\,(2t)^{l-2j},
\end{equation}
cf. \cite[Sec.~IX.3.1, formula (3)]{Vilenkin}.
A set of Gegenbauer polynomials $\bigl\{C_l^\lambda\bigr\}_{l\in\mathbb N_0}$ builds a complete orthogonal system on $[-1,1]$ with weight $(1-t^2)^{\lambda-1/2}$. Consequently, it is an orthogonal basis for zonal functions on the $2\lambda+1$--dimensional sphere.

Let $Q_l$ denote a polynomial on~$\mathbb{R}^{n+1}$ homogeneous of degree~$l$, i.e., such that $Q_l(az)=a^lQ_l(z)$ for all $a\in\mathbb R$ and $z\in\mathbb R^{n+1}$, and harmonic in~$\mathbb{R}^{n+1}$, i.e., satisfying $\nabla^2Q_l(z)=0$, then $Y_l(x)=Q_l(x)$, $x\in\mathcal S^n$, is called a hyperspherical harmonic of degree~$l$. The set of hyperspherical harmonics of degree~$l$ restricted to~$\mathcal S^n$ is denoted by $\mathcal H_l(\mathcal S^n)$. Hyperspherical harmonics of distinct degrees are orthogonal to each other. The number of linearly independent hyperspherical harmonics of degree~$l$ is equal to
$$
N=N(n,l)=\frac{(n+2l-1)(n+l-2)!}{(n-1)!\,l!}.
$$
Addition theorem states that
\begin{equation}\label{eq:addition_theorem}
C_l^\lambda(x\cdot y)=\frac{\lambda}{\lambda+l}\,\sum_{\kappa=1}^N\overline{Y_l^\kappa(x)}\,Y_l^\kappa(y)
\end{equation}
for any orthonormal set $\{Y_l^\kappa\}_{\kappa=1,2,\dots,N(n,l)}$ of hyperspherical harmonics of degree~$l$ on~$\mathcal S^n$. A proof can be found in~\cite{HTFII}, but note that hyperspherical harmonics are normalized such that
\begin{equation}\label{eq:HTF_Ykl_normalization}
\int_{\mathcal S^n}\left|Y_k^l(x)\right|^2d\sigma(x)=1.
\end{equation}
The authors of~\cite{LV} cite addition theorem from~\cite{HTFII} without changing the constant, and the same mistake is repeated in~\cite{EBCK09}, although in both articles one has~$\Sigma_n$ on the right--hand side of~\eqref{eq:HTF_Ykl_normalization}.

In this paper, we will be working with the orthogonal basis for~$\mathcal L^2(\mathcal S^n)=\overline{\bigoplus_{l=0}^\infty\mathcal H_l}$, consisting of hyperspherical harmonics given by
\begin{equation}\label{eq:spherical_harmonics}
Y_l^k(x)=A_l^k\prod_{\tau=1}^{n-1}C_{k_{\tau-1}-k_\tau}^{\frac{n-\tau}{2}+k_\tau}(\cos\theta_\tau)\sin^{k_\tau}\!\theta_\tau\cdot e^{\pm ik_{n-1}\varphi}
\end{equation}
with $l=k_0\geq k_1\geq\dots\geq k_{n-1}\geq0$, $k$ being a sequence $(k_1,\dots,\pm k_{n-1})$ of integer numbers, and normalization constants
\begin{equation}\label{eq:Alk}
A_l^k=\left(\frac{1}{\Gamma\left(\frac{n+1}{2}\right)}\prod_{\tau=1}^{n-1}\frac{2^{n-\tau+2k_\tau-2}\,(k_{\tau-1}-k_\tau)!\,(n-\tau+2k_{\tau-1})\,\Gamma^2(\frac{n-\tau}{2}+k_\tau)}{\sqrt\pi\,\Gamma(n-\tau+k_{\tau-1}+k_\tau)}\right)^{1/2},
\end{equation}
compare~\cite[Sec.~IX.3.6, formulae (4) and (5)]{Vilenkin}. The set of nonincreasing sequences~$k$ in $\mathbb N_0^{n-1}\times\mathbb Z$ with elements bounded by~$l$ will be denoted by $\mathcal M_{n-1}(l)$. For $k=(0,0,\dots,0)$ we have
\begin{equation}\label{eq:Al0}
A_l^0=\sqrt{\frac{(n-2)!\,l!\,(n+2l-1)}{(n+l-2)!\,(n-1)}},
\end{equation}
cf. \cite[Sec.~IX.3.6, formula (6)]{Vilenkin} (alternatively, use doubling formula \cite[8.335.1]{GR} for the computation).

Further, Funk--Hecke theorem states that for $f\in\mathcal L_\lambda^1\bigl([-1,1])$ and $Y_l\in\mathcal H_l(\mathcal S^n)$, $l\in\mathbb{N}_0$,
\begin{equation}\label{eq:FH}\begin{split}
\int_{\mathcal S^n}&Y_l(y)\,f(x\cdot y)\,d\sigma(y)\\
&=Y_l(x)\cdot\frac{(4\pi)^\lambda\,l!\,\Gamma(\lambda)}{(2\lambda+l-1)!}\int_{-1}^1 f(t)\,C_l^\lambda(t)\left(1-t^2\right)^{\lambda-1/2}dt.
\end{split}\end{equation}

Every $\mathcal{L}^1(\mathcal S^n)$--function~$f$ can be expanded into Laplace series of hyperspherical harmonics by
$$
S(f;x)\sim\sum_{l=0}^\infty Y_l(f;x),
$$
where $Y_l(f;x)$ is given by
\begin{equation*}\begin{split}
Y_l(f;x)&=\frac{\Gamma(\lambda)(\lambda+l)}{2\pi^{\lambda+1}}\int_{\mathcal S^n}f(y)\,C_l^\lambda(x\cdot y)\,d\sigma(y)\\
&=\frac{\lambda+l}{\lambda\Sigma_n}\int_{\mathcal S^n}f(y)\,C_l^\lambda(x\cdot y)\,d\sigma(y).
\end{split}\end{equation*}
For zonal functions we obtain by~\eqref{eq:FH} the representation
$$
Y_l(f;t)=\widehat f(l)\,C_l^\lambda(t),\qquad t=\cos\theta,
$$
with Gegenbauer coefficients
$$
\widehat f(l)=c(l,\lambda)\int_{-1}^1 f(t)\,C_l^\lambda(t)\left(1-t^2\right)^{\lambda-1/2}dt,
$$
where
\begin{equation*}
c(l,\lambda)=\frac{2^{2\lambda-1}\Gamma^2(\lambda)(\lambda+l)\Gamma(l+1)}{\pi\,\Gamma(2\lambda+l)}
   =\frac{\Gamma(\lambda)\Gamma(2\lambda)(\lambda+l)\Gamma(l+1)}{\sqrt\pi\,\Gamma(\lambda+\frac{1}{2})\Gamma(2\lambda+l)},
\end{equation*}
compare \cite[p. 207]{LV}.
The series
\begin{equation}\label{eq:Gegenbauer_expansion}
\sum_{l=0}^\infty\widehat f(l)\,C_l^\lambda(t)
\end{equation}
is called Gegenbauer expansion of~$f$.

For $f,h\in\mathcal L^1(\mathcal S^n)$, $h$ zonal, their convolution $f\ast h$ is defined by
\begin{equation*}
(f\ast h)(x)=\frac{1}{\Sigma_n}\int_{\mathcal S^n}f(y)\,h(x\cdot y)\,d\sigma(y).
\end{equation*}
With this notation we have
$$
Y_l(f;x)=\frac{\lambda+l}{\lambda}\,\bigl(f\ast C_l^\lambda\bigr)(x),
$$
hence, the function $\frac{\lambda+l}{\lambda}\,C_l^\lambda$ is the reproducing kernel for~$\mathcal H_l(\mathcal S^n)$, and Funk--Hecke formula can be written as
$$
Y_l\ast f=\frac{\lambda}{\lambda+l}\,\widehat f(l)\,Y_l.
$$

Further, any function $f\in\mathcal L^2(\mathcal S^n)$ has a unique representation as a mean--convergent series
\begin{equation}\label{eq:Fourier_series}
f(x)=\sum_{l=0}^\infty\sum_{k\in\mathcal M_{n-1}(l)} a_l^k\,Y_l^k(x),\qquad x\in\mathcal S^n,
\end{equation}
where
$$
a_l^k=a_l^k(f)=\frac{1}{\Sigma_n}\int_{\mathcal S^n}\overline{Y_l^k(x)}\,f(x)\,d\sigma(x)=\left<Y_l^k,f\right>,
$$
for proof cf.~\cite{Vilenkin}. In analogy to the two-dimensional case, we call $a_l^k$ the Fourier coefficients of the function~$f$. Convolution with a zonal function can be then written as
$$
f\ast g=\sum_{l=0}^\infty\sum_{k\in\mathcal M(_{n-1}l)} \frac{\lambda}{\lambda+l}\,a_l^k(f)\,\widehat g(l)\,Y_l^k
$$
and for zonal functions the following relation
\begin{equation}\label{eq:hatfl_vs_al0f}
\widehat f(l)=A_l^0\cdot a_l^0(f)
\end{equation}
between Fourier and Gegenbauer coefficients holds with~$A_l^0$ given by~\eqref{eq:Al0} (compare formulae \eqref{eq:spherical_harmonics}, \eqref{eq:Gegenbauer_expansion}, and~\eqref{eq:Fourier_series}).

Young's inequality for spherical functions has the following form:
$$
\|f\ast h\|_{\mathcal L^r}\leq\|f\|_{\mathcal L^p}\cdot\|h\|_{\mathcal L_\lambda^q},\qquad\frac{1}{r}=\frac{1}{q}+\frac{1}{p}-1\geq0,
$$
for $f\in\mathcal L^p(\mathcal S^n)$ and $h\in\mathcal L_\lambda^q\bigl([-1,1]\bigr)$.

The set of rotations of~$\mathbb R^{n+1}$ is denoted by $SO(n+1)$. It is isomorphic to the set of square matrices of degree $n+1$ with determinant~$1$. The $n$--dimensional sphere can be identified with the class of left cosets of $SO(n+1)$ mod $SO(n)$,
$$
\mathcal S^n=SO(n+1)/SO(n),
$$
cf.~\cite[Sec.~I.2]{Vilenkin}. Wigner polynomials $\{T_l^{km}\}$, $l\in\mathbb N_0$, $k,m\in\mathcal M_{n-1}(l)$, on $SO(n+1)$ are given by the regular representation of~$SO(n+1)$ in~$\mathcal L^2(\mathcal S^n)$:
\begin{equation}\label{eq:representationSO(n+1)}
Y_l^k(g^{-1}x)=\sum_{m\in\mathcal M_{n-1}(l)}T_l^{km}(g)\,Y_l^m(x),
\end{equation}
$x\in\mathcal S^n$, $g\in SO(n+1)$. Note that the order of indices we use is opposite to the one from~\cite[Sec.~IX.4.1, formula~(1)]{Vilenkin}. Invariant subspaces of this representation are $\mathcal H_l$, further, the matrix $(T_l^{\kappa\mu})_{\kappa,\mu=1}^{N(n,l)}$ is unitary. The normalized Haar measure on $SO(n+1)$ is denoted by~$d\nu$. The set
$$
\left\{\sqrt{N(n,l)}\,T_l^{km}:\,l\in\mathbb N_0,\,k,m\in\mathcal M_{n-1}(l)\right\}
$$
is orthonormal on $(SO(n+1),d\nu)$ \cite[Sec.~I.4.3, Thm.~1]{Vilenkin}. For $m=(0,0,\dots,0)$, Wigner polynomials~$T_l^{k0}$ are given by
\begin{equation}\label{eq:Tlk0}
T_l^{k0}(g)=\overline{T_l^{0k}(g^{-1})}=\sqrt{\frac{(n-1)!\,l!}{(n+2l-1)(n+l-2)!}}\,\overline{Y_l^k(g^{-1}\hat e)}
\end{equation}
compare \cite[Sec.~IX.4]{Vilenkin}.

Zonal product of arbitrary $\mathcal L^2(\mathcal S^n)$--functions $f$ and $h$ was introduced in~\cite{EBCK09} as
$$
(f\hat\ast h)(x\cdot y)=\int_{SO(n+1)}f(g^{-1}x)\,h(g^{-1}y)\,d\nu(g),\qquad x,y\in\mathcal S^n,
$$
and it has the representation
\begin{equation}\label{eq:zonal_product_Fc}
(f\hat\ast h)(x\cdot y)=\sum_{l=0}^\infty\sum_{k\in\mathcal M_{n-1}(l)}\frac{a_l^k(f)\,a_l^k(h)}{N(n,l)}\,\frac{\lambda+l}{\lambda}\,C_l^\lambda(x\cdot y).
\end{equation}
A difference in constant towards~\cite[formula (4.6)]{EBCK09} is caused by the fact that addition theorem is used in the proof, and the factor~$\Sigma_n$ must be corrected.

For further details on this topic we refer to the textbooks~\cite{Vilenkin} and~\cite{AH12}.

\subsection{Singular integrals and approximate identities}

Spherical singular integrals were introduced in~\cite{LV} by Berens \emph{et al.}, inspired by some previous papers concerning the theory of singular integrals on the real line~\cite{pB60}, unit circle~\cite{pB61,SW58}, in $k$-dimensional Euclidean space~\cite{CZ55} or on the $k$-dimensional torus~\cite{NP68}.

\begin{df}\label{def:singular_integral} Let  $\{\mathcal{K}_\rho\}_{\rho\in\mathbb{R}_+}\subseteq\mathcal{L}_\lambda^1\bigl([-1,1]\bigr)$ be a family of kernels such that
\begin{equation}\label{eq:hatK0}
\widehat{\mathcal K_\rho}(0)=c(0,\lambda)\int_{-1}^1\mathcal K_\rho(t)\bigl(1-t^2\bigr)^{\lambda-1/2}dt=1,
\end{equation}
Then
\begin{equation}\label{eq:Irho}
I_\rho(f)=f\ast\mathcal K_\rho
\end{equation}
is called a spherical singular integral. The family $\{\mathcal K_\rho\}$ is called the kernel of a singular integral.
\end{df}

\begin{bfseries}Remark. \end{bfseries} In~\cite{LV} the limit $\rho\to\infty$ is used. Here, we follow the convention $\rho\to0^+$ from~\cite{FW}, which is more convenient for the examples we discuss.

Approximate identities (without using this notion) were studied in~\cite{LV}, originally understood as singular integrals having an additional property that
\begin{equation}\label{eq:approximate_identity}
\lim_{\rho\to0^+}\|I_\rho f-f\|_{\mathcal{X}}=0.
\end{equation}

This definition is used e.g. in \cite{sB09,BE10KBW,BE10WS3,FW,FW-C,FGS-book,AMVA08}. However, condition~\eqref{eq:hatK0} is necessary neither for the approximation property~\eqref{eq:approximate_identity} nor for the definition of spherical wavelets. Moreover, it is not satisfied by many wavelet families. Therefore, similarly as in~\cite{EBCK09}, we define approximate identities as follows.

\begin{df} Let a family $\{\mathcal{K}_\rho\}_{\rho\in\mathbb{R}_+}$  of integrable zonal functions satisfying~\eqref{eq:approximate_identity} be given. Then, the family $\{\mathcal K_\rho\ast\}_{\rho\in\mathbb{R}_+}$ forms an approximate identity with kernel $\{\mathcal{K}_\rho\}_{\rho\in\mathbb{R}_+}$.
\end{df}

\begin{bfseries}Remark. \end{bfseries}Condition~\cite[(3.11)]{EBCK09} is satisfied by Young's inequality, therefore, it is not necessary in the definition.

Another useful characterization of approximate identities is given in the next theorem, cf. \cite[Theorem 3.8]{EBCK09}.

\begin{thm}\label{thm:approximate_identity}
Assume that the kernel $\{\mathcal{K}_\rho\}_{\rho\in\mathbb{R}_+}\subseteq\mathcal{L}_\lambda^1\bigl([-1,1]\bigr)$ is uniformly bounded in~$\mathcal L_\lambda^1$--norm, i.e.,
\begin{equation}\label{eq:L1_norm}
\int_{-1}^1|\mathcal K_\rho(t)|\bigl(1-t^2\bigr)^{\lambda-1/2}dt\leq\mathfrak c
\end{equation}
uniformly in $\rho\in\mathbb{R}_+$ for a positive constant~$\mathfrak c$. Then, the corresponding integral~$I_\rho$, defined by~\eqref{eq:Irho} is an approximate identity in~$\mathcal X(\mathcal S^n)$ if and only if
\begin{equation}\label{eq:limKrho}
\lim_{\rho\to0^+}\widehat{\mathcal K_\rho}(l)=\frac{\lambda+l}{\lambda}
\end{equation}
for all $l\in\mathbb N_0$.\end{thm}

Examples are provided by Abel--Poisson singular integral with the kernel
\begin{equation*}
A_\rho(t)=\frac{1-r^2}{(1-2rt+r^2)^\lambda}=\sum_{l=0}^\infty\frac{\lambda+l}{\lambda}\,r^l\,C_l^\lambda(t),\qquad r=e^{-\rho},
\end{equation*}
and by Gauss--Weierstrass singular integral with the kernel
\begin{equation*}
W_\rho(t)=\sum_{l=0}^\infty\frac{\lambda+l}{\lambda}\, e^{-l(l+2\lambda)\rho}\,C_l^\lambda(t).
\end{equation*}
Their properties are discussed in detail in~\cite{LV}.

\subsection{Dilations via stereographic projection}\label{subs:dilations}

Following~\cite{AVn}, we can define dilations of spherical zonal functions by inverse stereographic projection from the south pole. In the stereographic projection~$S$ a point~$x$ on the sphere is mapped onto a point~$\xi$ on the tangent space at the north pole being the intersection point of this tangent space and the straight line going through~$x$ and the south pole. We have
$$
|\xi|=2\tan\frac{\theta}{2}.
$$
Then, usual dilation is performed in the tangent space, and the image of~$\xi$ is projected back onto the sphere. Explicitly we have
\begin{equation}\label{eq:thetaa}
\tan\frac{\theta^a}{2}=a\tan\frac{\theta}{2},\qquad a\in\mathbb R_+.
\end{equation}
The ($\mathcal L^1$--norm preserving) dilation is given by
\begin{equation}\label{eq:L1dilation}\begin{split}
f^a(\theta^a,\theta_2,\dots,\theta_{n-1},\varphi)&=(D_af)(\theta,\theta_2,\dots,\theta_{n-1},\varphi)\\
&=\mu(a,\theta)\cdot f(\theta,\theta_2,\dots,\theta_{n-1},\varphi),
\end{split}\end{equation}
where the $\mu$--factor is the Radon--Nikodym derivative which ensures the covariance of the measure $d\sigma$ under dilation,
$$
\sin^{2\lambda}\theta\,d\theta=\mu(a,\theta)\sin^{2\lambda}\theta^a\,d\theta^a.
$$
Consequently, we have
\begin{equation}\label{eq:fathetadsigma}
f^a(\theta^a,\theta_2,\dots,\theta_{n-1},\varphi)\,\sin^{2\lambda}\theta^a\,d\theta^a=f(\theta,\theta_2,\dots,\theta_{n-1},\varphi)\,\sin^{2\lambda}\theta\,d\theta.
\end{equation}
The relation~\eqref{eq:thetaa} implies for $t=\cos\theta$ and $t^a=\cos\theta^a$
\begin{equation}\label{eq:changeofvariables}
t^a=\frac{(1-a^2)+(a^2+1)t}{(1-a^2)t+(a^2+1)}.
\end{equation}

It is proven in~\cite{ADJV02} (compare also~\cite[Section~9.3]{AMVA08}) and in~\cite{sB09} that in the case of a two-, respectively three-dimensional sphere, a dilation operator defined in this manner creates a kernel of an approximate identity. Here we show that it holds in any dimension.

\begin{thm}
Let $f\in\mathcal C\bigl([-1,1])$ be given. Then, the family $\{f^a\}_{a\in\mathbb{R}_+}$ is a kernel of an approximate identity in~$\mathcal X(\mathcal S^n)$.
\end{thm}

\begin{bfseries}Proof. \end{bfseries} Condition~\eqref{eq:L1_norm} is satisfied by~\eqref{eq:fathetadsigma} (applied to the function~$|f|$; note that~$\mu$ must be positive). In order to prove~\eqref{eq:limKrho}, consider
$$
\widehat{f^a}(l)=c(l,\lambda)\int_{-1}^1f^a(t^a)\,C_l^\lambda(t^a)\,\left(1-(t^a)^2\right)^{\lambda-1/2}\,dt^a.
$$
With a change of variables~\eqref{eq:changeofvariables} and by the measure--invariance we obtain
$$
\widehat{f^a}(l)=c(l,\lambda)\int_{-1}^1f(t)\,C_l^\lambda\left(\frac{(1-a^2)+(a^2+1)t}{(1-a^2)t+(a^2+1)}\right)\left(1-t^2\right)^{\lambda-1/2}\,dt.
$$
This expression is uniformly bounded by
\begin{align*}
c(l,\lambda)\,&\max_{t\in[-1,1]}|f(t)|\,\max_{t\in[-1,1]}\bigl|C_l^\lambda(t)\bigr|\,\int_{-1}^1\left(1-t^2\right)^{\lambda-1/2}\,dt\\
&\leq\text{const}\cdot\max_{t\in[-1,1]}|f(t)|
\end{align*}
and since
$$
\lim_{a\to0^+}C_l^\lambda\left(\frac{(1-a^2)+(a^2+1)t}{(1-a^2)t+(a^2+1)}\right)=C_l^\lambda(1)=\binom{2\lambda+l-1}{l},
$$
compare~\cite[formula~8.937.4]{GR}, we have
$$
\lim_{a\to0^+}\widehat{f^a}(l)
   =\binom{2\lambda+l-1}{l}\cdot c(l,\lambda)\underbrace{\int_{-1}^1f(t)\left(1-t^2\right)^{\lambda-1/2}\,dt}_{=\widehat{f}(0)/c(0,\lambda)}=\frac{\lambda+l}{\lambda}
$$
and hence, \eqref{eq:limKrho} is satisfied for all $l\in\mathbb N_0$.\hfill$\Box$

\section{Continuous wavelet transform -- bilinear theory}\label{sec:bilinearwv}

We first recall the definition of bilinear wavelets, given in~\cite{EBCK09}, with a slight modification concerning conjugation, such that also complex families of functions can be considered.

\begin{df}\label{def:bilinear_wavelets} Let $\alpha:\mathbb R_+\to\mathbb R_+$ be a weight function. A family $\{\Psi_\rho\}_{\rho\in\mathbb R_+}\subseteq\mathcal L^2(\mathcal S^n)$ is called bilinear spherical wavelet if it satisfies the following admissibility conditions:
\begin{enumerate}
\item for $l\in\mathbb{N}_0$
\begin{equation}\label{eq:admbwv1}
\sum_{\kappa=1}^{N(n,l)}\int_0^\infty \bigl|a_l^\kappa(\Psi_\rho)\bigr|^2\,\alpha(\rho)\,d\rho=N(n,l),
\end{equation}
\item for $R\in\mathbb{R}_+$ and $x\in\mathcal S^n$
\begin{equation}\label{eq:admbwv2}
\int_{\mathcal S^n}\left|\int_R^\infty(\overline{\Psi_\rho}\hat\ast\Psi_\rho)(x\cdot y)\,\alpha(\rho)\,d\rho\right|d\sigma(y)\leq\mathfrak c
\end{equation}
with $\mathfrak c$ independent of~$R$.
\end{enumerate}
\end{df}

Again, the factor~$\Sigma_n$ must be corrected (with respect to the formulae from~\cite{EBCK09}), because of the constant correction in addition theorem.

\begin{df}\label{def:bilinear_wt} Let $\{\Psi_\rho\}_{\rho\in\mathbb R_+}$ be a spherical wavelet. Then, the spherical wavelet transform
$$
\mathcal W_\Psi\colon\mathcal L^2(\mathcal S^n)\to\mathcal L^2(\mathbb R_+\times SO(n+1))
$$
is defined by
$$
\mathcal W_\Psi f(\rho,g)=\frac{1}{\Sigma_n}\int_{\mathcal S^n}\overline{\Psi_\rho(g^{-1}x)}\,f(x)\,d\sigma(x).
$$
\end{df}

As it is stated in \cite[Theorem~5.3]{EBCK09}, the wavelet transform defined in this manner is invertible (in $\mathcal L^2$--sense) by
\begin{equation}\label{eq:bwt_synthesis}
f(x)=\int_0^\infty\!\!\int_{SO(n+1)}\mathcal W_\Psi f(\rho,g)\,\Psi_\rho(g^{-1}x)\,d\nu(g)\,\alpha(\rho)\,d\rho.
\end{equation}

This reconstruction formula yields
$$
f=\int_0^\infty\left[f\ast\left(\overline{\Psi_\rho}\hat\ast\Psi_\rho\right)\right]\,\alpha(\rho)\,d\rho,
$$
i.e., the function $\overline{\Psi_\rho}\hat\ast\Psi_\rho$ is the reproducing kernel of the image $\mathcal W_\Psi(\mathcal S^n)$ of the wavelet transform.

\begin{bfseries}Remark 1. \end{bfseries}$\rho$ is the index of the wavelet family. In this general setting it is not the parameter of stereographic dilation.

\begin{bfseries}Remark 2. \end{bfseries}Note that the zero-mean condition
$$
a_0^0(\Psi_\rho)=0
$$
is not required. This is an important difference to the classical wavelet theory on~$\mathbb R^n$ or to the definitions from~\cite{AV} and~\cite{AVn}. However, in order to justify the usage of the word \emph{wavelet}, one could add this requirement, similarly as in~\cite{FW} and~\cite{FGS-book}. Moreover, zero-mean wavelets are most common in computer applications. One can also consider wavelets of order~$m\in\mathbb N_0$, i.e., such with $m+1$ vanishing moments
$$
a_l^k(\Psi_\rho)=0\qquad\text{for }l=0,\,1,\,\dots,m,\,k\in\mathcal M_{n-1}(l),
$$
and~\eqref{eq:admbwv1} holding for $l=m+1,\,m+2,\,\dots$ (cf. \cite[Sec.~10.2]{FGS-book}). In this case, the reconstruction formula is valid for functions~$f$ with $m+1$ vanishing moments. For the proof, we take the function $\sum_{l=0}^m\frac{\lambda+l}{\lambda}\,C_l^\lambda+\Xi_R$ as the kernel of an approximate identity.

\begin{bfseries}Remark 3.\end{bfseries}\label{rem:L2_boundedness} The $\mathcal L^2$--boundedness of the wavelet is necessary for the existence of hyperspherical harmonics representation, $SO(n+1)$--rotation and for the zonal product to be defined. In the case of zonal wavelets, the wavelet transform reduces to a convolution over~$\mathcal S^n$, cf. formula~\eqref{eq:zonal_bilinear_wt}, and both the wavelet transform and its inverse converge for $\mathcal L^1$--wavelets.

Reconstruction formula is proven in the original paper~\cite{EBCK09}, here, we show that the wavelet transform is an isometry.

\begin{thm}\label{thm:isometry} Let $\{\Psi_\rho\}_{\rho\in\mathbb R_+}$ be a spherical wavelet and $f,h\in\mathcal L^2(\mathcal S^n)$. Then,
\begin{equation}\label{eq:isometry}
\left<\mathcal W_\Psi f,\mathcal W_\Psi h\right>=\left<f,h\right>
\end{equation}
where the scalar product in the wavelet phase space is given by
$$
\left<F,H\right>_{\mathcal L^2(\mathbb R_+\times SO(n+1))}=\int_0^\infty\int_{SO(n+1)}\overline{F(\rho,g)}\,H(\rho,g)\,d\nu(g)\,\alpha(\rho)\,d\rho.
$$
\end{thm}

\begin{bfseries}Proof.\end{bfseries}
\begin{align*}
&\left<\mathcal W_\Psi f,\mathcal W_\Psi h\right>=\int_0^\infty\int_{SO(n+1)}\overline{\mathcal W_\Psi f(\rho,g)}\,\mathcal W_\Psi h(\rho,g)\,d\nu(g)\,\alpha(\rho)\,d\rho\\
&=\int_0^\infty\!\!\int_{SO(n+1)}\int_{\mathcal S^n}\!\Psi(g^{-1}x)\,\overline{f(x)}d\sigma(x)\,
   \int_{\mathcal S^n}\!\overline{\Psi(g^{-1}y)}\,h(y)\,d\sigma(y)\,d\nu(g)\,\alpha(\rho)\,d\rho.
\end{align*}
A change of integration order yields
\begin{align*}
&\left<\mathcal W_\Psi f,\mathcal W_\Psi h\right>\\
&=\int_{\mathcal S^n}\!\int_{\mathcal S^n}\!\int_0^\infty\!\!\!
   \underbrace{\int_{SO(n+1)}\Psi_\rho(g^{-1}x)\,\overline{\Psi_\rho(g^{-1}y)}\,d\nu(g)}_{\Psi_\rho\hat\ast\overline{\Psi_\rho}(x\cdot y)}
   \alpha(\rho)\,d\rho\,\overline{f(x)}\,d\sigma(x)\,g(y)\,d\sigma(y).
\end{align*}
Now, we replace the zonal product of $\Psi_\rho$'s by its Gegenbauer expansion~\eqref{eq:zonal_product_Fc} and obtain
\begin{align*}
&\left<\mathcal W_\Psi f,\mathcal W_\Psi h\right>\\
&=\int_{\mathcal S^n}\!\int_{\mathcal S^n}\!\int_0^\infty
   \frac{\sum_{l,k}\bigl|a_l^k(\Psi_\rho)\bigr|^2}{N(n,l)}\,\frac{\lambda+l}{\lambda}\,C_l^\lambda(x\cdot y)\,\alpha(\rho)\,d\rho\,\overline{f(x)}\,d\sigma(x)\,g(y)\,d\sigma(y).
\end{align*}
According to the property \eqref{eq:admbwv1}, the inner integral is equal to $\frac{\lambda+l}{\lambda}\,C_l^\lambda(x\cdot y)$, i.e., it is the reproducing kernel of the space of harmonic polynomials of degree~$l$. Consequently, the integral with respect to~$x$ yields~$\overline{f(y)}$, and finally, we obtain~\eqref{eq:isometry}.\nolinebreak\hfill $\Box$

\subsection{Euclidean limit}

A very important feature of wavelets is their scaling behavior. It is expected that for small scales they resemble wavelets over Euclidean space. Actually, this property is one of the starting points for the probably most popular spherical wavelet definition of Antoine and Vandergheynst~\cite{AVn}. Also Holschneider~\cite{mH96}, although not so exacting in other requirements (see the widely criticized \emph{ad hoc} introduction of the scale parameter or existence of the inverse transform \emph{whenever the integral makes sense} \cite[Sec.~2.2.2]{HCM03}), finds the so--called Euclidean limit property very important. On the other hand, small scale behavior of the wavelet itself has not been studied so far by the authors whose spherical wavelet definitions are based on spherical singular integrals \cite{sB09,BE10KBW,BE10WS3,sE11,FW,FW-C,FGS-book}. In this Subsection we fill this gap by showing that wavelets given by
$$
\Psi_\rho=\sum_{l=0}^\infty\sum_{k\in\mathcal M_{n-1}(l)}a_l^k(\Psi_\rho)\,Y_l^k
$$
under some mild conditions on the coefficients~$a_l^k$ behave for small scales like over the Euclidean space. The investigation and the proof are inspired by~\cite{mH96}.

\begin{thm}\label{thm:euclidean_limit}Let a spherical wavelet $\Psi_\rho\subseteq\mathcal L^2(\mathcal S^n)$ with
\begin{equation}\label{eq:alk_cond_Eucl_lim}
a_l^k(\Psi_\rho)=\frac{1}{l^{k_1-1}A_l^k}\,\mathcal O\left(\psi_k(l\rho)\right),\qquad\rho\to0,
\end{equation}
$l\in\mathbb N_0$, $k=(k_1,k_2,\dots,k_{n-1})\in\mathcal M_{n-1}(l)$, $k_1\leq K$, be given, with  $\psi_k\in\mathcal L^2(\mathbb R_+,\,t^{n-1}\,dt)$  -- a piecewise smooth function satisfying
\begin{equation}\label{eq:smallforsmallscales}
\rho^n\sum_{l=0}^{[c/\rho]}l^{n-1}\,\psi_k(l\rho)<\epsilon,
\end{equation}
for some $c>0$, $\epsilon\ll1$, and $\rho<\rho_0$. Further, let
$$
\lim_{\rho\to0}a_l^k(\Psi_\rho)=0
$$
for $k_1\geq K$. Then there exists a  square integrable function $F:\,\mathbb R^n\to\mathbb C$ such that
$$
\lim_{\rho\to0}\,\rho^n\,\Psi_\rho\left(S^{-1}(\rho\xi)\right)=F(\xi)
$$
holds point-wise for every $\xi\in\mathbb R^n$. $S^{-1}$ denotes the inverse stereographic projection.
\end{thm}

\begin{bfseries}Proof.\end{bfseries} Let $\xi\in\mathbb R^n$ be given in spherical coordinates, i.e.,
$$
\xi=(r,\theta_2,\theta_3,\dots,\theta_{n-1},\varphi)
$$
with $r=|\xi|$ and $(\theta_2,\theta_3,\dots,\theta_{n-1},\varphi)$ -- spherical coordinates of an $n-1$--dimensional sphere of radius~$r$. Further, let~$\theta$ denote the $\theta_1$--coordinate of $x:=S^{-1}(\rho\xi)$. Then
$$
\theta=2\arctan\frac{\rho r}{2}=\rho r+\mathcal O\left(\rho^3\right),\qquad\rho\to0,
$$
and hence,
\begin{equation}\label{eq:theta_asymptotics}
\cos\theta=\cos(\rho r)+\mathcal O(\rho^4),\qquad\rho\to0.
\end{equation}
The other coordinates remain unaffected by the (inverse) stereographic projection, i.e., $\theta_j(x)=\theta_j(\rho\xi)$ for $j=2,3,\dots,n-1$, and $\varphi(x)=\varphi(\rho\xi)$. Now, hyperspherical harmonics are functions of the form
\begin{equation*}\begin{split}
Y_l^k(\theta_1,\dots,\theta_{n-1},\varphi)=&\,C_{l-k_1}^{\lambda+k_1}(\cos\theta_1)\sin^{k_1}\!\theta_1
\cdot\frac{A_l^k}{A_{k_1}^{\widetilde k}}Y_{k_1}^{\widetilde k}(\theta_2,\dots,\theta_{n-1},\varphi),
\end{split}\end{equation*}
with $\widetilde k=(k_2,\dots,k_{n-1})$ and $Y_{k_1}^{\widetilde k}$ -- hyperspherical harmonics over the sphere~$\mathcal S^{n-1}$, compare~\eqref{eq:spherical_harmonics}. Further, from~\eqref{eq:C_jako_szereg} and~\eqref{eq:theta_asymptotics} it follows that
$$
C_{l-k_1}^{\lambda+k_1}(\cos\theta)=C_{l-k_1}^{\lambda+k_1}(\cos(\rho r))+\text{\scriptsize$\mathcal O$}(\rho^3),\qquad\rho\to0.
$$
On the other hand,
$$
\sin\theta=\sin(\rho r)+\mathcal O(\rho^3),\qquad\rho\to0,
$$
and, hence, we have the following representation for $\rho^n\Psi_\rho$:
\begin{equation}\label{eq:rhonPsin}\begin{split}
\rho^n&\Psi_\rho(\theta,\theta_2,\dots,\theta_{n-1},\varphi)\\
&=\rho^n\sum_{l=0}^\infty\sum_{k\in\mathcal M_{n-1}(l)}a_l^k(\Psi_\rho)\,[Y_l^k(\rho r,\theta_2,\dots,\theta_{n-1},\varphi)+\mathcal O(\rho^3)],\qquad\rho\to0.
\end{split}\end{equation}
Further, for any $l\in\mathbb N_0$ and $\xi\in\mathcal S^n$
$$
\sum_{k\in\mathcal M_{n-1}(l)}|Y_l^k(\xi)|^2=\frac{\lambda+l}{\lambda}\,C_l^\lambda(1)=\frac{\lambda+l}{\lambda}\,\binom{2\lambda+l-1}{l},
$$
and thus,
$$
\sum_{k\in\mathcal M_{n-1}(l)}|Y_l^k(\xi)|^2=\mathcal O(l^{2\lambda}),\qquad l\to\infty.
$$
compare~\eqref{eq:addition_theorem} and~\cite[formula~8.937.4]{GR}. Therefore, we obtain from~\eqref{eq:rhonPsin} by Schwarz inequality
$$
\left|\rho^n\Psi_n\left(S^{-1}(\rho r)\right)\right|\leq\rho^n\sum_{l=0}^\infty\sqrt{\sum_{k\in\mathcal M_{n-1}(l)}|a_l^k(\Psi_\rho)|^2}
   \cdot\text{const}\cdot  l^\lambda\left(1+\mathcal O(\rho^3)\right)],\quad\rho\to0.
$$
It follows from~\eqref{eq:smallforsmallscales} that the $l$'s contributing to the sum in~\eqref{eq:rhonPsin} in the limit $\rho\to0$  get large since the terms with $l\rho<c$ can be neglected.

For further considerations, denote by~$\widetilde k$ the sequence $(k_2,\dots,k_{n-1})$ and note that
$$
\sum_{l=0}^\infty\sum_{k\in\mathcal M_{n-1}(l)}=\sum_{l=0}^\infty\sum_{k_1=0}^{\min\{l,K\}}\sum_{\widetilde k\in\mathcal M_{n-2}(k_1)}.
$$
A change of summation order yields
\begin{equation*}\begin{split}
\rho^n&\Psi_\rho(\theta,\theta_2,\dots,\theta_{n-1},\varphi)=\rho^n\sum_{k_1=0}^K\sum_{l=k_1}^\infty\sum_{\widetilde k\in\mathcal M_{n-2}(k_1)} a_l^k(\Psi_\rho)\\
&\cdot\left[C_{l-k_1}^{\lambda+k_1}\left(\cos(\rho r)\right)\sin^{k_1}(\rho r)
   \cdot\frac{A_l^k}{A_{k_1}^{\widetilde k}}\,Y_{k_1}^{\widetilde k}(\theta_2,\dots,\theta_{n-1},\varphi)+\mathcal O(\rho^3)\right],\quad\rho\to0,
\end{split}\end{equation*}
supposed the series converges absolutely. For the Gegenbauer polynomials we now use the representation
$$
C_l^\lambda(t)=\frac{\Gamma(l+2\lambda)\,\Gamma(\lambda+\frac{1}{2})}{\Gamma(2\lambda)\,\Gamma(\lambda+l+\frac{1}{2})}\,\,P_l^{(\lambda-\frac{1}{2},\lambda-\frac{1}{2})}(t)
$$
with $ P_l^{(\lambda,\lambda)}$ -- Jacobi polynomials, cf.~\cite[formula 8.962]{GR}. Thus, we have
\begin{equation}\label{eq:rhonPsirho}\begin{split}
\rho^n&\Psi_\rho(\theta,\theta_2,\dots,\theta_{n-1},\varphi)=\rho^n\sum_{k_1=0}^K\sum_{l=k_1}^\infty\sum_{\widetilde k\in\mathcal M_{n-2}(k_1)} a_l^k(\Psi_\rho)\\
&\cdot\Biggl[\frac{\Gamma(l+k_1+2\lambda)\Gamma(k_1+\lambda+\frac{1}{2})}{\Gamma(2k_1+2\lambda)\Gamma(\lambda+l+\frac{1}{2})}
   P_{l-k_1}^{(k_1+\lambda-\frac{1}{2}, k_1+\lambda-\frac{1}{2})}\left(\cos(\rho r)\right)\sin^{k_1}(\rho r)\\
&\cdot\frac{A_l^k}{A_{k_1}^{\widetilde k}}\,Y_{k_1}^{\widetilde k}(\theta_2,\dots,\theta_{n-1},\varphi)+\mathcal O(\rho^3)\Biggr],\qquad\rho\to0.
\end{split}\end{equation}
Since
\begin{equation}\label{eq:limit_Gamma_Gamma}
\frac{\Gamma(l+k_1+2\lambda)}{\Gamma(\lambda+l+\frac{1}{2})}=\mathcal O\left(l^{k_1+\lambda-\frac{1}{2}}\right),\qquad l\to\infty,
\end{equation}
and
$$
\lim_{l\to\infty}l^{-\lambda} P_l^{(\lambda,\lambda)}\left(\cos\frac{r}{l}\right)=\left(\frac{r}{2}\right)^{-\lambda}J_\lambda(r),
$$
where~$J_\lambda$ is the $\lambda^{\text{th}}$ order Bessel function, compare~\cite[formula 8.722]{GR}, the right-hand-side of~\eqref{eq:rhonPsirho} behaves in limit $\rho\to0$ like
\begin{align}
\sum_{k_1=0}^K&\frac{C_{k_1}\Gamma(k_1+\lambda+\frac{1}{2})\,\rho^{2\lambda+1}}{\lambda\,\Gamma(2k_1+2\lambda)}
   \sum_{l=k_1}^\infty\frac{\psi_k(l\rho)}{l^{k_1-1}}\cdot l^{2k_1+2\lambda-1}\left(\frac{l\rho r}{2}\right)^{\frac{1}{2}-k_1-\lambda}\notag\\
&\cdot J_{k_1+\lambda-\frac{1}{2}}(l\rho r)\cdot(\rho r)^{k_1}
   \sum_{\widetilde k\in\mathcal M_{n-2}(k_1)}\frac{1}{A_{k_1}^{\widetilde k}}\,Y_{k_1}^{\widetilde k}(\theta_2,\dots,\theta_{n-1},\varphi)\notag\\
&=\sum_{k_1=0}^K\frac{2^{k_1+\lambda-\frac{1}{2}}C_{k_1}\,\Gamma(k_1+\lambda+\frac{1}{2})}{\lambda\,\Gamma(2k_1+2\lambda)\,r^{\lambda-\frac{1}{2}}}
   \sum_{l=k_1}^\infty\rho\,(l\rho)^{\lambda+\frac{1}{2}}\psi_k(l\rho)\,J_{k_1+\lambda-\frac{1}{2}}(l\rho r)\notag\\
&\cdot\sum_{\widetilde k\in\mathcal M_{n-2}(k_1)}\frac{Y_{k_1}^{\widetilde k}(\theta_2,\dots,\theta_{n-1},\varphi)}{A_{k_1}^{\widetilde k}},\label{eq:seriesrhoto0}
\end{align}
where $C_{k_1}$ is the proportionality constant in~\eqref{eq:limit_Gamma_Gamma}.
Now, the function $t\mapsto t^{\lambda+1/2}\psi_k(t)$ is in $\mathcal L^1$ and piecewise continuous, and $|J_{\lambda-\frac{1}{2}}(\cdot)|\leq1$. Therefore, the series in~\eqref{eq:seriesrhoto0} converges absolutely to the Lebesgue integral. By substituting $t$ for $l\rho$ and $dt$ for $\rho$, we can write
\begin{equation*}\begin{split}
\lim_{\rho\to0}\rho^n&\Psi_\rho(\theta,\theta_2,\dots,\theta_{n-1},\varphi)\\
&=\sum_{k_1=0}^K\frac{2^{k_1+\lambda-\frac{1}{2}}C_{k_1}\Gamma(k_1+\lambda+\frac{1}{2})}{\lambda\,\Gamma(2k_1+2\lambda)\,r^{\lambda-\frac{1}{2}}}
   \int_0^\infty t^{\lambda+\frac{1}{2}}\,\psi_k(t)\,J_{k_1+\lambda-\frac{1}{2}}(tr)\,dt\notag\\
&\cdot\sum_{\widetilde k\in\mathcal M_{n-2}(k_1)}\frac{Y_{k_1}^{\widetilde k}(\theta_2,\dots,\theta_{n-1},\varphi)}{A_{k_1}^{\widetilde k}}
\end{split}\end{equation*}
According to~\cite[Theorem~IV.3.3]{SW}, the expression
$$
\frac{2\pi}{r^{k_1+\lambda-\frac{1}{2}}}\int_0^\infty t^{k_1+\lambda+\frac{1}{2}}\,\frac{\psi_k(t)}{t^{k_1}}\,J_{k_1+\lambda-\frac{1}{2}}(tr)\,dt
$$
is equal to the Fourier transform of the function $F_k:\,\mathbb R^{n+2k_1}\to\mathbb C,\,\xi\mapsto\frac{\psi_k(|\xi|)}{|\xi|^{k_1}}$. Consequently, we have
\begin{equation}\label{eq:limes_rhoPsi}\begin{split}
\lim_{\rho\to0}\rho^n\Psi_\rho(\theta,\theta_2,\dots,\theta_{n-1},\varphi)
   =&\sum_{k_1=0}^K\frac{2^{k_1+\lambda-\frac{1}{2}}C_{k_1}\Gamma(k_1+\lambda+\frac{1}{2})\,r^{k_1}\widehat{F_k}(r)}{2\pi\lambda\,\Gamma(2k_1+2\lambda)}\\
&\cdot\sum_{\widetilde k\in\mathcal M_{n-2}(k_1)}\frac{Y_{k_1}^{\widetilde k}(\theta_2,\dots,\theta_{n-1},\varphi)}{A_{k_1}^{\widetilde k}}.
\end{split}\end{equation}
Since the sum over $\widetilde k\in\mathcal M_{n-2}(k_1)$ is bounded by a constant, square integrability of $r^{k_1}\widehat{F_k}(r)$ for each~$k_1$ ensures that the limit of $\rho^n\Psi_\rho$ belongs to~$\mathcal L^2(\mathbb R^n)$. Note that although $\widehat{F_k}$ was computed as a function over $\mathbb R^{n+2k_1}$, it is a radial one, and in~\eqref{eq:limes_rhoPsi} it occurs as a function of $r=|\xi|$, $\xi\in\mathbb R^n$. Its $\mathcal L^2$--norm is equal to
\begin{align*}
\|r^{k_1}\widehat{F_k}\|_{\mathcal L^2(\mathbb R^n)}^2&=\int_0^\infty|\widehat{F_k}|^2\,r^{n+2k_1-1}\,dr
   =\|\widehat{F_k}\|_{\mathcal L^2(\mathbb R^{n+2k_1})}\\
&=C\,\|F_k\|_{\mathcal L^2(\mathbb R^{n+2k_1})}=C\int_0^\infty\left|\frac{\psi_k(r)}{r^{k_1}}\right|^2\,r^{n+2k_1-1}\,dr\\
&=C\,\int_0^\infty|\psi_k(r)|^2\,r^{n-1}\,dr,
\end{align*}
and it is finite by assumption. Therefore, $\lim_{\rho\to0}\rho^n\Psi_\rho\in\mathcal L^2(\mathbb R^n)$.\hfill$\Box$

\begin{bfseries}Example. \end{bfseries}Fourier coefficients of directional Poisson wavelets~\cite{HH09} over $\mathcal S^2$ are equal to
$$
a_l^k(g_\rho^d)=\rho^de^{-\rho l}(-1)^k\sqrt{\frac{2l+1}{4\pi}}\,\sum_{j=0}^d\left(\binom{d}{j}\left(\frac{l}{2}\right)^d(-1)^j+\mathcal O(l^{d-1})\right)\,\delta_{k,2j-d},
$$
cf. \cite[p. 07512-6]{HH09}. Since the function $\gamma_d\colon t\mapsto t^de^{-t}$ is bounded, the second summand in this expression
$$
\rho\cdot\mathcal O\left(\gamma_{d-1}(\rho l)\right)
$$
vanishes for $\rho\to0$. We obtain
\begin{align*}
a_l^k(g_\rho^d)&=(-1)^k\sqrt{\frac{2l+1}{4\pi}}\cdot\mathcal O\!\left((\rho l)^de^{-\rho l}\,\sum_{j=0}^d\binom{d}{j}(-1)^j\delta_{k,2j-d}\right)\\
&=(-1)^{k+\frac{d+k}{2}}\binom{d}{\frac{d+k}{2}}\sqrt{\frac{2l+1}{4\pi}}\,\mathcal O\left(\gamma_d(\rho l)\right),\qquad\rho\to0,
\end{align*}
for $d$, $k$ with the same parity and such that $|k|\leq d$, and $a_l^k(g_\rho^d)=\mathcal O(\rho)$, $\rho\to0$, in other cases.

Now, since
$$
A_l^k=2^k\Gamma\left(k+\frac{1}{2}\right)\sqrt{\frac{(2l+1)(l-k)!}{(l+k)!}}=\mathcal O(l^{-k+\frac{1}{2}}),\qquad l\to\infty,
$$
see formula~\eqref{eq:Alk}, and $k=k_1$, it is straightforward to verify that the coefficients $a_l^k(g_\rho^d)$ satisfy condition~\eqref{eq:alk_cond_Eucl_lim} for $\psi_k=(-1)^{k+\frac{d+k}{2}}\binom{d}{\frac{d+k}{2}}\gamma_d$ (when  $2|(d-k)$), respectively $\psi_k=0$ (otherwise), with $K=d$. Thus, the theorem applies and Euclidean limit of directional Poisson wavelets exists as an $\mathcal L^2$--function. It is computed explicitly in~\cite[Sec.~IV.B]{HH09}, compare also~\cite[Thm.~6.1]{IIN14DW} and the Remark following it.

\subsection{Zonal wavelets}

In the case of zonal wavelets, Definition~\ref{def:bilinear_wavelets} and Definition~\ref{def:bilinear_wt} reduce to the following ones (by formula~\eqref{eq:hatfl_vs_al0f}). This task is the subject of~\cite[Section~5.5]{EBCK09}, but note that formula~(5.2) is false (Fourier coefficient of a function is supposed to be dependent on the argument).

\begin{df}\label{def:zonal_bilinear_wavelet}
A subfamily $\{\Psi_\rho\}_{\rho\in\mathbb{R}_+}$ of the space $\mathcal{L}_\lambda^1\bigl([-1,1]\bigr)$ is called zonal spherical wavelet if it satisfies the following admissibility conditions:
\begin{enumerate}
\item for $l\in\mathbb{N}_0$
\begin{equation}\label{eq:admwv1}
\int_0^\infty\bigl|\widehat{\Psi_\rho}(l)\bigr|^2\,\alpha(\rho)\,d\rho=\left(\frac{\lambda+l}{\lambda}\right)^2,
\end{equation}
\item for $R\in\mathbb{R}_+$
\begin{equation}\label{eq:admwv2}
\int_{-1}^1\left|\int_R^\infty\bigl(\overline{\Psi_\rho}\ast\Psi_\rho\bigr)(t)\,\alpha(\rho)\,d\rho\right|\,\left(1-t^2\right)^{\lambda-1/2}dt\leq\mathfrak c
\end{equation}
with $\mathfrak c$ independent of~$R$.
\end{enumerate}
\end{df}

\begin{bfseries}Remark 2. \end{bfseries}If $\alpha(\rho)=\frac{1}{\rho}$, the case studied in~\cite{FW}, \cite{FW-C} or in~\cite{FGS-book}, \eqref{eq:admwv1} implies
$$
\lim_{\rho\to0^+}\widehat\Psi_\rho(l)=0
$$
for $l\in\mathbb{N}_0$.

\begin{bfseries}Remark 3.\end{bfseries} If there exists a function~$\psi$ such that $\widehat\Psi_\rho(l)=\frac{\lambda+l}{\lambda}\,\psi(l\rho)$, then the condition~\eqref{eq:admwv1} can be written as
\begin{equation*}
\int_0^\infty\bigl|\psi(t)\bigr|^2\,\frac{dt}{t}=1.
\end{equation*}
In this case we call~$\psi$ the generating function of the wavelet. Note that by~\eqref{eq:hatfl_vs_al0f}, condition~\eqref{eq:alk_cond_Eucl_lim} is then satisfied, and the wavelet has the Euclidean limit property if~\eqref{eq:smallforsmallscales} holds.

\begin{df}
The $x$--rotation operator $R_x$, $x\in\mathcal{S}^n$, acting on the space of zonal functions, is given by
$$
R_x:\,\Psi\mapsto\Psi\bigl(\left<x,\cdot\right>\bigr).
$$
For a wavelet family $\{\Psi_\rho\}_{\rho\in\mathbb{R}_+}$, the $\rho$--dilation operator~$D_\rho$ is given by
\begin{equation}\label{eq:dilation}
D_\rho:\,\Psi_1\mapsto\Psi_\rho.
\end{equation}
By $\Psi_{\rho,x}$ we denote the dilated and rotated wavelet
$$
\Psi_{\rho,x}:\,y\mapsto\Psi_{\rho,x}(y)=\Psi_\rho(x\cdot y)=R_x D_\rho\Psi(y),\qquad y\in\mathcal{S}^n.
$$
Further, the inner product on $\mathbb{R}_+\times\mathcal{S}^n$ is given by
$$
\left<f,h\right>_{\mathcal{L}^2\left(\mathbb{R}_+\times\mathcal{S}^n\right)}
   =\frac{1}{\Sigma_n}\int_0^\infty\!\!\int_{\mathcal{S}^n}\,\overline{f(\rho,x)}\,h(\rho,x)\,d\sigma(y)\,\alpha(\rho)\,d\rho
$$
and the space $\mathcal{L}^2\left(\mathbb{R}_+\times\mathcal{S}^n\right)$ consists of functions~$f$ such that
$$
\|f\|_{\mathcal{L}^2\left(\mathbb{R}_+\times\mathcal{S}^n\right)}^2=\left<f,f\right>_{\mathcal{L}^2\left(\mathbb{R}_+\times\mathcal{S}^n\right)}<\infty.
$$
\end{df}

\begin{bfseries}Remark. \end{bfseries}Rotation and dilation are independent of each other -- in contrast to the $\mathbb R^n$--case and to the stereographic dilation. In fact, dilation operator~\eqref{eq:dilation} can be applied to a family of functions, not to a single one. However, a family $\{\Psi_\rho\}$ can be generated by a single function, for instance by the stereographic dilation.

\begin{df}\label{def:swt}
Let $\{\Psi_\rho\}_{\rho\in\mathbb{R}_+}$ be a spherical wavelet. Then, the spherical wavelet transform
$$
\mathcal{W}_\Psi:\,\mathcal{L}^2\left(\mathcal{S}^n\right)\to\mathcal{L}^2\left(\mathbb{R}_+\times\mathcal{S}^n\right)
$$
is defined by
\begin{equation}\label{eq:zonal_bilinear_wt}
\mathcal{W}_\Psi f(\rho,x)=\left(f\ast\overline{\Psi_\rho}\right)(x).
\end{equation}
\end{df}

By Young's inequality, the $\mathcal L^1$--boundedness of the wavelet is sufficient for the convergence of~\eqref{eq:zonal_bilinear_wt}, however, if the wavelet is an $\mathcal L^2$--function, the wavelet transform can be written as
$$
\mathcal{W}_\Psi f(\rho,x)=\left<\Psi_{\rho,x},f\right>.
$$
The wavelet transform is an isometry and it is invertible in $\mathcal L^2$--sense by
\begin{equation}\label{eq:reconstruction}
f(x)=\frac{1}{\Sigma_n}\int_0^\infty\!\!\!\int_{\mathcal S^n}\mathcal W_\Psi f(\rho,y)\,\Psi_{\rho,y}(x)\,d\sigma(y)\,\alpha(\rho)\,d\rho
\end{equation}
(the constant in the reconstruction formula given for zonal wavelets in~\cite[Section~5.2]{EBCK09} is valid for wavelet transform over~$\mathcal S^3$). It is a consequence of isometry and invertibility of the wavelet transform with respect to nonzonal wavelets, but both properties can be also proven directly (compare proofs in two-- and three--dimensional cases \cite{FW,FW-C,FGS-book,sB09,BE10KBW}). For the reconstruction formula one introduces the function
\begin{equation*}
\Xi_R=\int_R^\infty\bigl(\overline{\Psi_\rho}\ast\Psi_\rho\bigr)\,\alpha(\rho)\,d\rho
\end{equation*}
and shows that it is a kernel of an approximate identity.

\label{dis:L1_boundedness}Again, by Young's inequality $\mathcal L^1$--boundedness is sufficient for the convergence of this integral, therefore, in contrast to the nonzonal case, $\mathcal L^2$--boundedness is not required, compare Remark~3 on page~\pageref{rem:L2_boundedness}. For $\mathcal L^2$--wavelets the reconstruction formula can be written in the form:
$$
f(x)=\int_0^\infty\left<\overline{\Psi_{\rho,x}},\mathcal W_\Psi f\right>\alpha(\rho)\,d\rho.
$$

\begin{bfseries}Remark 1.\end{bfseries} The reproducing kernel of $\mathcal W_\Psi(\mathcal S^n)$ is $\overline{\Psi_\rho}\ast\Psi_\rho$.

\begin{bfseries}Remark 2.\end{bfseries} Condition
\begin{equation}\label{eq:add_cond}
\sum_{l=0}^\infty\frac{\lambda+l}{\lambda}\int_R^\infty\bigl|\widehat{\Psi_\rho}(l)\bigr|^2\,\alpha(\rho)\,d\rho<\infty,
\end{equation}
analogous to~(2.3) in~\cite{FW}, is necessary neither for the existence of the wavelet transform, nor for the reconstruction formula. The existence of~$\Xi_R$ as an $\mathcal L^1$--function is ensured by condition~\eqref{eq:admwv2}.

\subsection{Spherical scaling function and wavelets corresponding to an approximate identity}\label{subs:sc_fct}

In the case of zonal wavelets, a scaling function can be defined in a similar way as in~\cite{FW}.

\begin{df} Let $\{\Psi_\rho\}_{\rho\in\mathbb{R}_+}$ be a spherical wavelet. Then, the corresponding spherical scaling function $\{\Phi_R\}_{R\in\mathbb{R}_+}$ is defined by
$$
\Phi_R(t)=\sum_{l=0}^\infty\widehat{\Phi_R}(l)\,C_l^\lambda(t),\qquad t\in[-1,1],
$$
where
$$
\widehat{\Phi_R}(l)=\left(\int_R^\infty\bigl|\widehat{\Psi_\rho}(l)\bigr|^2\,\alpha(\rho)\,d\rho\right)^{1/2}
$$
for $l\in\mathbb{N}_0$.
\end{df}

\begin{bfseries}Remark 1.\end{bfseries} Note that
$$
\Phi_R\ast\Phi_R=\Xi_R.
$$

\begin{bfseries}Remark 2. \end{bfseries}For $\Psi_\rho$ satisfying~\eqref{eq:add_cond}, we have $\Phi_R\in\mathcal L^2(\mathcal S^n)$.

Conversely, to a suitable kernel of an approximate identity, a spherical wavelet can be associated, cf. \cite[Theorem~6.1]{EBCK09}.

\begin{thm}\label{thm:zonal_wv_from_kernel}
Let a kernel $\{\Phi_R\}_{R\in\mathbb R_+}$ of a uniformly bounded approximate identity be given with Gegenbauer coefficients which are differentiable with respect to~$R$ and monotonically decreasing in~$R$. Moreover, assume that
$$
\lim_{R\to\infty}\widehat{\Phi_R}(l)=0
$$
for $l\in\mathbb N$. Then, the associated spherical wavelet $\{\Psi_\rho\}_{\rho\in\mathbb R_+}$ is given by
$$
\widehat{\Psi_\rho}(l)=\left(-\frac{1}{\alpha(\rho)}\,\frac{d}{d\rho}\,\bigl|\widehat{\Phi_\rho}(l)\bigr|^2\right)^{1/2}
$$
for $l\in\mathbb N_0$, $\rho\in\mathbb R_+$.
\end{thm}

\begin{bfseries}Remark.\end{bfseries} Uniform boundedness of a family of functions is necessary for~\eqref{eq:limKrho} to be a defining property of an approximate identity, cf.~Theorem~\ref{thm:approximate_identity}, and in the proof of \cite[Theorem~6.1]{EBCK09} this case is assumed to hold. Note, however, that an approximate identity does not need to be uniformly bounded, compare Definition~\ref{def:singular_integral} and formula~\eqref{eq:approximate_identity}.

Spherical wavelets, associated to singular integrals are Abel--Poisson wavelet with Gegenbauer coefficients of the kernel
$$
\widehat{\Psi_\rho}(l)=\frac{\lambda+l}{\lambda}\,\sqrt{2l\rho}\,e^{-l\rho},\qquad l\in\mathbb N_0,
$$
and Gauss--Weierstrass wavelet given by
$$
\widehat{\Psi_\rho}(l)=\frac{\lambda+l}{\lambda}\,\sqrt{2l(l+2\lambda)\rho}\,e^{-l(l+2\lambda)\rho},\qquad l\in\mathbb N_0,
$$
compare also~\cite{FW}. Note that in both examples the wavelets are of order~$0$.

A much more interesting task is a construction of nonzonal wavelets from an approximate identity, described in~\cite{EBCK09} (we have corrected the factor~$\Sigma_n$ both in the definition of admissible vectors and in the theorem in order to obtain compatibility with the corrected addition theorem).

\begin{df}\label{def:admissible_weight_vectors}
Let $w_l\in\mathbb R^{N(n,l)}$ be a vector with components $w_l(\kappa)$, $\kappa=1,2,\dots,N(n,l)$. Is is said to be an admissible weight vector if it satisfies
$$
\sum_{\kappa=1}^{N(n,l)}|w_l(\kappa)|^2=N(n,l).
$$
\end{df}

\begin{thm}\label{thm:bsw_from_AI} Let $\{w_l\}_{l\in\mathbb N}$ be a family of admissible weight vectors, and let $\{\Phi_R\}_{R\in\mathbb R_+}$ be a kernel of an approximate identity satisfying assumptions of Theorem~\ref{thm:zonal_wv_from_kernel}. The family of functions given by
$$
\Psi_\rho(x)=\sum_{l=0}^\infty\sum_{\kappa=0}^{N(n,l)}\left(-\frac{1}{\alpha(\rho)}\,\frac{d}{d\rho}\,\bigl|\widehat{\Phi_\rho}(l)\bigr|^2\right)^{\frac{1}{2}}
   \frac{\lambda}{\lambda+l}\,w_l(\kappa)\,Y_l^\kappa(x)
$$
is a bilinear spherical wavelet.
\end{thm}

\begin{bfseries}Remark. \end{bfseries}Weight vectors~$w_l$, $l\in\mathbb N$ with
$$
w_l(\kappa)=\sqrt{N(n,l)}\,\delta_{\kappa,0}
$$
yield zonal vectors.

We would like to point out that admissible weight vectors could also be dependent on the scale parameter~$\rho$. The following enhancement of Definition~\ref{def:admissible_weight_vectors} does not influence the validity of Theorem~\ref{thm:bsw_from_AI} (with $w_l$ replaced by~$\omega_l(\rho)$).

\begin{df}
Let $\omega_l(\rho)\in\mathbb R^{N(n,l)}$ be a vector with components $w_l(\rho,\kappa)$, $\kappa=1,2,\dots,N(n,l)$, $\rho\in\mathbb R_+$. It is said to be an admissible weight vector if it satisfies
$$
\sum_{\kappa=1}^{N(n,l)}|\omega_l(\rho,\kappa)|^2=N(n,l)\qquad\text{for any }\rho\in\mathbb R_+.
$$
\end{df}

Note, however, that using scale--dependent weight vectors makes it more difficult to construct wavelets with Euclidean limit property.

\section{Linear wavelet transform}\label{sec:linearwv}

Linear wavelet transform based on singular integrals was introduced for two-dimensional spheres in~\cite{FW}, compare also~\cite{FW-C} and the references therein. A generalization to three dimensions can be found in~\cite{sB09} and in \cite{BE10KBW}. In these papers, only zonal wavelets are considered. Nonzonal wavelets are studied in \cite{EBCK09} an in \cite{BE10KBW} only in the bilinear version. Here, we would like to introduce linear wavelet transform with respect to nonzonal wavelets.

\begin{df}\label{def:linear_wavelets} Let $\alpha:\mathbb R_+\to\mathbb R_+$ be a weight function. A family $\{\Psi_\rho^L\}_{\rho\in\mathbb R_+}\subseteq\mathcal L^2(\mathcal S^n)$ is called linear spherical wavelet if it satisfies the following admissibility conditions:
\begin{enumerate}
\item for $l\in\mathbb{N}_0$
\begin{equation}\label{eq:admlwv1}
A_l^0\cdot\int_0^\infty\,a_l^0(\Psi_\rho^L)\,\alpha(\rho)\,d\rho=\frac{\lambda+l}{\lambda},
\end{equation}
\item for $R\in\mathbb{R}_+$ and $x\in\mathcal S^n$
\begin{equation}\label{eq:admlwv2}
\int_{\mathcal S^n}\left|\int_R^\infty\Psi_\rho^L(x\cdot y)\,\alpha(\rho)\,d\rho\right|d\sigma(y)\leq\mathfrak c
\end{equation}
with $\mathfrak c$ independent of~$R$.
\end{enumerate}
\end{df}

\begin{df}\label{def:linear_wt} Let $\{\Psi_\rho^L\}_{\rho\in\mathbb R_+}$ be a linear spherical wavelet. Then, the linear spherical wavelet transform
$$
\mathcal W_\Psi^L\colon\mathcal L^2(\mathcal S^n)\to\mathcal L^2\left(\mathbb R_+\times SO(n+1)\right)
$$
is defined by
\begin{equation}\label{eq:lwt}
\mathcal W_\Psi^L f(\rho,g)=\frac{1}{\Sigma_n}\int_{\mathcal S^n}\Psi_\rho^L(g^{-1}x)\,f(x)\,d\sigma(x).
\end{equation}
\end{df}

\begin{thm}(Reconstruction formula) Let $\{\Psi_\rho^L\}_{\rho\in\mathbb R_+}$ be a linear wavelet and $f\in\mathcal L^2(\mathcal S^n)$. Then
\begin{equation}\label{eq:linear_rec_formula}
f(x)=\int_0^\infty\!\!\int_{SO(n)}\mathcal W_\Psi^L f(\rho,\hat x\tilde g)\,d\nu(\tilde g)\,\alpha(\rho)\,d\rho.
\end{equation}
in $\mathcal L^2$--sense, where $\hat x$ denotes any fixed element of~$SO(n+1)$ satisfying $\hat x\hat e=x$.
\end{thm}

\begin{bfseries}Proof. \end{bfseries}Substitute~\eqref{eq:lwt} in the inner integral in~\eqref{eq:linear_rec_formula} and change the order of integration,
\begin{align*}
\int_{SO(n)}\mathcal W_\Psi^L &f(\rho,\hat x\tilde g)\,d\nu(\tilde g)
   =\frac{1}{\Sigma_n}\int_{SO(n)}\int_{\mathcal S^n}\Psi_\rho^L\left((\hat x\tilde g)^{-1}y\right)\,f(y)\,d\sigma(y)\,d\nu(\tilde g)\\
&=\frac{1}{\Sigma_n}\int_{\mathcal S^n}\int_{SO(n)}\Psi_\rho^L(\tilde g^{-1}\hat x^{-1}y)\,d\nu(\tilde g)f(y)\,d\sigma(y).
\end{align*}
We show that the function
$$
\Upsilon_{\!\rho}(z)\colon=\int_{SO(n)}\Psi_\rho^L(\tilde g^{-1}z)\,d\nu(\tilde g),\quad z\in\mathcal S^n,
$$
is zonal. Let~$T$ be the regular representation of $SO(n)\subseteq SO(n+1)$ in~$\mathcal L^2(\mathcal S^n)$, and \mbox{$\tilde g_0\in SO(n)$}. Then
$$
T(\tilde g_0)\Upsilon_{\!\rho}(z)=\int_{SO(n)}\Psi_\rho^L\left(\tilde g_0^{-1}\tilde g^{-1}z\right)\,d\nu(\tilde g)=\int_{SO(n)}\Psi_\rho^L\left((\tilde g\tilde g_0)^{-1}z\right)\,d\nu(\tilde g).
$$
By measure invariance we thus obtain
$$
T(\tilde g_0)\Upsilon_{\!\rho}(z)=\int_{SO(n)}\Psi_\rho^L\left((\tilde g\tilde g_0)^{-1}z\right)\,d\nu(\tilde g\tilde g_0)=\Upsilon_{\!\rho}(z).
$$
Fourier coefficients~$a_l^k$ of this function vanish for $k\ne(0,0,\dots,0)$, and for $k=(0,0,\dots,0)$ they are equal to
\begin{align*}
a_l^0(\Upsilon_{\!\rho})&=\left<Y_l^0,\Upsilon_{\!\rho}\right>
   =\left<Y_l^0,\int_{SO(n)}\sum_{l=0}^\infty\sum_{k\in\mathcal M_{n-1}(l)}a_l^k(\Psi_\rho^L)Y_l^k(\tilde g^{-1}\circ)\,d\nu(\tilde g)\right>.
\end{align*}
By \eqref{eq:representationSO(n+1)} we have
\begin{align*}
a_l^0(\Upsilon_{\!\rho})&=\sum_{k\in\mathcal M_{n-1}(l)}a_l^k(\Psi_\rho^L)\int_{SO(n)}\left<Y_l^0,\sum_{m\in\mathcal M_{n-1}(l)}T_l^{k m}(\tilde g)Y_l^ m\right>d\nu(\tilde g)\\
&=\sum_{k\in\mathcal M_{n-1}(l)}a_l^k(\Psi_\rho^L)\int_{SO(n)}T_l^{k0}(\tilde g)\,d\nu(\tilde g),
\end{align*}
and further, \eqref{eq:Tlk0} yields
\begin{equation*}
a_l^0(\Upsilon_{\!\rho})=\sum_{k\in\mathcal M_{n-1}(l)}a_l^k(\Psi_\rho^L)\sqrt{\frac{(n-1)!\,l!}{(n+2l-1)(n+l-2)!}}\,\int_{SO(n)}\overline{Y_l^k(\tilde g^{-1}\hat e)}\,d\nu(\tilde g).
\end{equation*}
Since any $\tilde g\in SO(n)$ leaves~$\hat e$ invariant, $\nu\left(SO(n)\right)=1$, and
$$
Y_l^k(\hat e)=\begin{cases}A_l^k\,C_l^\lambda(1)&\text{for }k=(0,0,\dots,0),\\0&\text{otherwise},\end{cases}
$$
compare~\eqref{eq:spherical_harmonics}, we obtain by \cite[formula 8.937.4]{GR}
\begin{align*}
a_l^0(\Upsilon_{\!\rho})=&a_l^0(\Psi_\rho^L)\sqrt{\frac{(n-1)!\,l!}{(n+2l-1)(n+l-2)!}}\\
&\cdot\sqrt{\frac{(n-2)!\,l!\,(n+2l-1)}{(n+l-2)!\,(n-1)}}\,\binom{n+l-2}{l}=a_l^0(\Psi_\rho^L).
\end{align*}

Now, choose a positive number~$R$. Then
$$
\int_R^\infty\!\!\int_{SO(n)}\mathcal W_\Psi^L f(\rho,\hat x\tilde g)\,d\nu(\tilde g)\,\alpha(\rho)\,d\rho=\bigl(f\ast\Phi_R^L\bigr)(x)
$$
for $\Phi_R^L$ with
$$
\widehat{\Phi_R^L}(t)=\int_R^\infty\widehat{\Upsilon_{\!\rho}}(l)\,\alpha(\rho)\,d\rho.
$$
It is now straightforward to show that the function $\Phi_R^L$ is a kernel of an approximate identity, and hence the theorem holds.\hfill$\Box$

\begin{bfseries}Remark 1. \end{bfseries}Similarly as in~\cite{sB09} and in the case of bilinear wavelets, the zero-mean condition is not required, compare also discussion in~\cite[p.~231]{FGS-book}. Again, one can also consider wavelets of order~$m$, $m\in\mathbb N_0$, i.e., with $m+1$ vanishing moments and~\eqref{eq:admlwv1} holding for $l=m+1,\,m+2,\,\dots$. In this case, the reconstruction formula is valid for functions~$f$ with $m+1$ vanishing moments. For the proof, we take the function $\sum_{l=0}^m\frac{\lambda+l}{\lambda}\,C_l^\lambda+\Phi_R^L$ as the kernel of an approximate identity.\\
\begin{bfseries}Remark 2. \end{bfseries}$\mathcal L^1$--boundedness of the wavelet is sufficient for the convergence of both the wavelet transform and its inverse, compare Remark~3 on page~\pageref{rem:L2_boundedness} and discussion for bilinear wavelets on page~\pageref{dis:L1_boundedness}.

\begin{cor}The reproducing kernel of the linear wavelet transform is the function
$$
\Upsilon_{\!\rho}=\int_{SO(n)}\Psi_\rho^L(\tilde g^{-1}\,\circ)\,d\nu(\tilde g).
$$
\end{cor}
\begin{bfseries}Proof.\end{bfseries} A substitution of~\eqref{eq:lwt} into~\eqref{eq:linear_rec_formula} yields
$$
f(x)=\int_0^\infty\,\frac{1}{\Sigma_n}\int_{\mathcal S^n}f(y)\cdot\Upsilon_{\!\rho}(\hat x^{-1}y)\,d\sigma(y)\,\alpha(\rho)\,d\rho.
$$
Since~$\Upsilon_{\!\rho}$ is zonal and $\hat x^{-1}x=\hat e$, we obtain $\Upsilon_{\!\rho}(\hat x^{-1}y)=\Upsilon_{\!\rho}(x\cdot y)$ and hence
$$
f=\int_0^\infty(f\ast\Upsilon_{\!\rho})\,\alpha(\rho)\,d\rho.
$$
\hfill$\Box$

Linear wavelet transform is in general not an isometry, compare the proof of Theorem~\ref{thm:isometry}, but the wavelets may posses the Euclidean limit property.

\subsection{Zonal linear wavelets}

In the case of zonal wavelets Definition~\ref{def:linear_wavelets} and Definition~\ref{def:linear_wt} reduce to the following ones.

\begin{df}\label{def:zonal_lsw} Let $\{\Psi_\rho^L\}_{\rho\in\mathbb R_+}$ be a subfamily of~$\mathcal L_\lambda^1\bigl([-1,1]\bigr)$ such that the following admissibility conditions are satisfied:
\begin{enumerate}
\item for $l\in\mathbb{N}_0$
$$
\int_0^\infty\widehat{\Psi_\rho^L}(l)\,\alpha(\rho)\,d\rho=\frac{\lambda+l}{\lambda},
$$
\item for $R\in\mathbb{R}_+$
$$
\int_{-1}^1\left|\int_R^\infty\Psi_\rho^L(t)\,\alpha(\rho)\,d\rho\right|\,\left(1-t^2\right)^{\lambda-1/2}dt\leq\mathfrak c
$$
with $\mathfrak c$ independent from~$R$.
\end{enumerate}
Then $\{\Psi_\rho^L\}_{\rho\in\mathbb R_+}$ is called a spherical linear wavelet. The associated wavelet transform
$$
\mathcal W_\Psi^L\colon\mathcal L^2(\mathcal S^n)\to\mathcal L^2(\mathbb R_+\times\mathcal S^n)
$$
is defined by
$$
\mathcal{W}_\Psi^L f(\rho,x)=\bigl(f\ast\overline{\Psi_\rho^L}\bigr)(x).
$$
\end{df}

In the case of an $\mathcal L^2$--wavelet, the wavelet transform can be written as
$$
\mathcal{W}_\Psi^L f(\rho,x)=\left<\Psi_{\rho,x}^L,f\right>_{\mathcal{L}^2(\mathcal{S}^n)}.
$$

In the reconstruction formula no integral over $SO(n)$ is needed.

\begin{thm}\label{thm:zonal_lwv_reconstruction} Let $\{\Psi_\rho^L\}_{\rho\in\mathbb R_+}$ be a linear wavelet and $f\in\mathcal L^2(\mathcal S^n)$. Then
$$
f(x)=\frac{1}{\Sigma_n}\int_0^\infty\mathcal{W}_\Psi^L f(\rho,x)\,\alpha(\rho)\,d\rho
$$
in $\mathcal L^2$--sense.
\end{thm}

Consequently, the reproducing kernel of~$\mathcal W_\Psi^L(\mathcal S^n)$ is equal to~$\Psi_\rho^L$.

\subsection{Spherical scaling function and zonal linear wavelets corresponding to an approximate identity}

Similarly as in the bilinear case, a scaling function can be associated to a wavelet.

\begin{df} Let $\{\Psi_\rho^L\}_{\rho\in\mathbb{R}_+}$ be a spherical zonal linear wavelet. Then, the corresponding spherical scaling function $\{\Phi_R\}_{R\in\mathbb{R}_+}$ is defined by
$$
\Phi_R^L(t)=\sum_{l=0}^\infty\widehat{\Phi_R^L}(l)\,C_l^\lambda(t),\quad t\in[-1,1],
$$
where
$$
\widehat{\Phi_R^L}(l)=\int_R^\infty\widehat{\Psi_\rho^L}(l)\,\alpha(\rho)\,d\rho
$$
for $l\in\mathbb{N}_0$.
\end{df}

Consequently, the scaling function itself is a kernel of an approximate identity.

Similarly as in Section~\ref{subs:sc_fct}, we can construct linear wavelets from kernels of approximate identities.

\begin{thm}
Let a kernel $\{\Phi_R\}_{R\in\mathbb R_+}$ of a uniformly bounded approximate identity be given with Gegenbauer coefficients which are differentiable with respect to~$R$. Moreover, assume that
$$
\lim_{R\to\infty}\widehat{\Phi_R}(l)=0
$$
for $l\in\mathbb N_0$. Then, the associated spherical linear wavelet $\{\Psi_\rho^L\}_{\rho\in\mathbb R_+}$ is given by
$$
\widehat{\Psi_\rho^L}(l)=-\frac{1}{\alpha(\rho)}\,\frac{d}{d\rho}\,\widehat{\Phi_\rho}(l)
$$
for $l\in\mathbb N_0$, $\rho\in\mathbb R_+$.
\end{thm}

\begin{bfseries}Proof. \end{bfseries}Analogous to the proof of \cite[Theorem~6.1]{EBCK09}.

\begin{bfseries}Remark 1. \end{bfseries}Contrary to the bilinear case, it is not necessary that Gegenbauer coefficients of the approximate identity are monotonously decreasing, as it is required in~\cite[Corollary~10.1.7]{FGS-book}. This property is not requested in~\cite[Definition~2.7]{BE10WS3}.

\begin{bfseries}Remark 2. \end{bfseries}Condition
$$
\sum_{l=0}^\infty\frac{\lambda+l}{\lambda}\int_R^\infty\widehat{\Psi_\rho^L}(l)\,\alpha(\rho)\,d\rho<\infty,
$$
for $R\in\mathbb R_+$, analogous to condition~(iii) in Definition~7.1 from~\cite{FW}, ensures that positivity of Gegenbauer coefficients of the kernel is equivalent to its positive definiteness, compare~\cite[Theorem~5.7.4]{FGS-book}, but it is necessary neither for the definition of the wavelets nor for the reconstruction formula.

Abel--Poisson spherical linear wavelet, derived from Abel--Poisson kernel, given by
$$
\widehat{\Psi_\rho^L}(l)=\frac{\lambda+l}{\lambda}\,l\rho\,e^{-l\rho},
$$
is a multipole Poisson wavelet of order~$d=1$, cf. e.g.~\cite{HI07}. It is shown in~\cite{IIN14PW} that also multipole wavelets of higher order are linear wavelets in the sense of Definition~\ref{def:zonal_lsw}.
Gauss--Weierstrass spherical linear wavelet, derived from Gauss--Weierstrass kernel, has Gegenbauer coefficients
$$
\widehat{\Psi_\rho^L}(l)=\frac{\lambda+l}{\lambda}\,l(l+2\lambda)\rho\,e^{-l(l+2\lambda)\rho}.
$$
In all these cases we deal with wavelets of order~$0$.

\subsection{Construction of nonzonal wavelets}
There is a big freedom in a construction of nonzonal linear wavelets. This is a consequence of condition~\eqref{eq:admlwv1} which is a constraint only on coefficients~$a_l^0$ of a wavelet family. Other coefficients must be chosen such that condition~\eqref{eq:admlwv2} is satisfied. Thus, additional requirements, e.g., on $\mathcal L^2$--norm or angular selectivity (cf. \cite[Sec.~9.2.4.2]{AMVA08}), can be requested.

\section{Relationship to other $n$--dimensional isotropic spherical wavelets}\label{sec:other_wavelets}\label{sec:other}

One of the main goals of this paper is to show that many of the existing zonal wavelet constructions can be treated as bilinear wavelets derived from approximate identities.

\subsection{Holschneider's wavelets}
One of the first attempts to define continuous wavelet transform for the $2$--sphere was done by Holschneider in~\cite{mH96}. The wavelet transform is performed in the same way as in the present paper, and the wavelets are admissible functions according to the following definition.

\begin{df}\label{df:Holschneider}A family of functions $\{\Psi_\rho\}_{\rho\in\mathbb R_+}\subseteq\mathcal C^\infty(\mathcal S^2)$ is admissible if it satisfies the conditions:\\
a) There exists a constant $c_\Psi$ such that for all $f\in\mathcal L^2(\mathcal S^2)$ with $\int_{\mathcal S^2}f=0$ we have
$$
\int_{\mathbb R_+\times SO(3)}|\mathcal W_\Psi f(\rho,g)|^2\,\frac{d\nu(g)\,d\rho}{\rho}=c_\Psi\int_{\mathcal S^2}|f(x)|^2\,d\sigma(x).
$$
b) For all $\alpha>0$ there exists a finite constant $c_\alpha$ such that for $\rho>1$ we have
$$
|\mathcal W_\Psi f(\rho,g)|\leq c_\alpha\rho^{-\alpha}.
$$
c) Euclidean limit property holds.
\end{df}
The wavelet synthesis is in principle performed by~\eqref{eq:bwt_synthesis}, however, with respect to another wavelet family~$\{\Omega_\rho\}$. Admissibility conditions for $\{\Psi_\rho\}$ and $\{\Omega_\rho\}$ to be an analysis--reconstruction pair are boundedness of the wavelet transform, boundedness of the wavelet synthesis and inversion formula, i.e., a requirement that~\eqref{eq:bwt_synthesis} holds.

There are two weak points of this construction that we want to discuss. The first one is the \emph{ad hoc} choice of the scale parameter, an accusation that has been made by many authors, see e.g. \cite{AMVA08,AV,AVn,CFK06,mF09,MHML}. The other one is the fact that wavelets are not defined intrinsically. In the case of zonal wavelets, it is explicitly stated that the wavelet synthesis is performed by (in principle)~\eqref{eq:reconstruction} \emph{whenever this integral makes sense}, cf. \cite[Sec.~2.2.2]{HCM03}.

It is worth noting that the only spherical wavelets that have been used in practise by Holschneider are Poisson multipole wavelets \cite{HCM03,CPMDHJ05,HI07,IH10} given by
$$
\widehat{\Psi_\rho^d}(l)=C(l)\,\frac{\lambda+l}{\lambda}\,(\rho l)^d\,e^{-\rho l},\qquad d\in\mathbb{N}.
$$
(Note that Abel--Poisson wavelet can be thought of as a Poisson wavelet with $d=\frac{1}{2}$.) It is shown in~\cite{IIN14PW} that these wavelets are indeed wavelets derived from a kernel of an approximate identity (both linear and bilinear). Similarly, directional Poisson wavelets~\cite{HH09} are bilinear wavelets derived from an approximate identity, however, a different reconstruction family is required~\cite{IIN14DW}.

Therefore, it is possible to replace Definition~\ref{df:Holschneider} by definitions from Section~\ref{sec:bilinearwv} (Holschneider's wavelets are wavelets with a generating function). This movement retrieves both of the above described drawbacks.

\subsection{Antoine and Vandergheynst' wavelets}

The most investigated spherical wavelet construction is the one derived purely group--theoretical~\cite{AVn}. A whole spectrum of papers was published on two--dimensional wavelets (by authors such as Antoine, Bogdanova, Ferreira, Jacques, McEwen, Rosca, Vandergheynst, Wiaux and others), and also $n$--dimensional wavelets are quite popular \cite{CFK06,mF08,mF09}.

In this approach, wavelets are obtained from an admissible mother wavelet via dilations on the tangent space given by~\eqref{eq:L1dilation} or
$$
f^a(\theta^a,\theta_2,\dots,\theta_{n-1},\varphi)=(D^af)(\theta)=\sqrt{\mu(a,\theta)}f(\theta,\theta_2,\dots,\theta_{n-1},\varphi)
$$
for an $\mathcal L^2$--norm preserving dilation.

It is shown in~\cite[Section~5.3]{EBCK09} that bilinear wavelets from Definition~\ref{def:bilinear_wavelets} satisfy admissibility conditions from~\cite{AVn}:
\begin{equation}\label{eq:AVcondition}
\int_{\mathbb R_+\times SO(n+1)}\left|\left<f,U(\sigma(\rho,g))\Psi_1\right>\right|^2\,d\nu(g)\,\alpha(\rho)\,d\rho<\infty,
\end{equation}
where the representation of the corresponding section over $\mathbb R_+\times SO(n+1)$ is defined by
\begin{align*}
U(\sigma(\rho,g))\Psi_1&=(T(g)\circ D(\rho))\Psi_1,\\
D(\rho)\Psi_1&=\Psi_\rho,\\
T(g)\Psi_\rho&=\Psi_\rho(g^{-1}\circ).
\end{align*}
This statement does not mean that bilinear wavelets derived from an approximate identity can be treated as a special case of Antoine and Vandergheynst' wavelets. The difference lies in the way how dilation is defined.

It is easy to see from the proof of~\eqref{eq:AVcondition} that linear wavelets from Definition~\ref{def:linear_wavelets} in general do not fulfill this condition.

On the other hand, it can be asked if group--theoretically derived wavelets can satisfy conditions of Definition~\ref{def:bilinear_wavelets} or Definition~\ref{def:linear_wavelets}. In the case of zonal wavelets, one of the possibilities is to verify whether for a given wavelet its scaling function (defined in the same way as in the present paper) is a uniformly bounded approximate identity (with monotonously decreasing Gegenbauer coefficients for the bilinear case). Note however, that Theorem~\ref{thm:approximate_identity} cannot be used for checking it (even in the case of an $\mathcal L^1$--norm preserving dilation) since the scaling behavior of wavelets~\eqref{eq:L1dilation} is not inherited by their scaling function.

\subsection{Mexican needlets}

Recently, a new construction has been proposed by Geller and Mayeli in~\cite{GM09a,GM09b}. The wavelets are kernels of the convolution operator $f(t\Delta^{\!\ast})$, where $0\ne f\in\mathcal S(\mathbb R_+)$, $f(0)\ne0$, and~$\Delta^{\!\ast}$ denotes the Laplace--Beltrami operator on a manifold. In the case of the sphere this leads to zonal wavelets of the form
\begin{equation}\label{eq:Mneedlets}
K_\rho(\hat e,y)=\frac{1}{\Sigma_n}\sum_{l=0}^\infty f\left(\rho^2l(l+2\lambda)\right)\frac{\lambda+l}{\lambda}\,C_l^\lambda(y)
\end{equation}
In~\cite{aM10} the name \emph{Mexican needlets} is introduced for kernels of the form~\eqref{eq:Mneedlets} with $f_r(s)=s^re^{-s}$, $r\in\mathbb N$. Similarly as in the case of Poisson wavelets~\cite[Sec.~8]{IIN14PW} one can easily prove that normalized Mexican needlets
$$
K_\rho^r(y)=2^r\,\sqrt{\frac{2}{\Gamma(2r)}}\,\sum_{l=0}^\infty\left(\rho^2l(l+2\lambda)\right)^r e^{-\rho^2l(l+2\lambda)}\,\frac{\lambda+l}{\lambda}\,C_l^\lambda(y)
$$
satisfy~\eqref{eq:admwv1} for $\alpha(\rho)=\frac{1}{\rho}$, the measure used in~\cite{GM09a,GM09b}. Condition~\eqref{eq:admwv2} can be proven along the same lines as boundedness of the kernel in~\cite[Proof of Thm.~5.7]{IIN14DW}. Consequently, Mexican needlets are also a special case of zonal bilinear wavelets derived from an approximate identity. Further, although Geller and Mayeli propose only a bilinear construction, it can be shown that
$$
\widetilde K_\rho^r(y)=\frac{2}{\Gamma(r)}\,\sum_{l=0}^\infty\left(\rho^2l(l+2\lambda)\right)^r e^{-\rho^2l(l+2\lambda)}\,\frac{\lambda+l}{\lambda}\,C_l^\lambda(y)
$$ are zonal linear wavelets according to the theory presented here.

\subsection{Ebert's diffusive wavelets}

Diffusive wavelets introduced in Ebert's PhD--thesis~\cite{sE11} are wavelets derived from an approximate identity and satisfying some additional conditions such as positivity and a semigroup condition. Therefore, zonal spherical diffusive wavelets are a special case of wavelets described in the present paper. Examples that Ebert gives are Gauss--Weierstrass and Abel--Poisson wavelets. It is worth noting in this place that Poisson multipole wavelets are not diffusive.

\end{document}